\documentstyle[12pt]{article}
\def\proof{\noindent{\bf Proof:}\hskip10pt}
\def\QED{\hfill $\Box$ \smallskip}

\font\tenmath=msbm10 scaled 1200
\font\sevenmath=msbm7 scaled 1200
\font\fivemath=msbm5 scaled 1200
\newfam\mathfam \textfont\mathfam=\tenmath
\scriptfont\mathfam=\sevenmath \scriptscriptfont\mathfam=\fivemath
\def\math{\fam\mathfam}
\begin{document}
\def \\ { \cr }
\def\R{{\math R}}
\def\N{{\math N}}
\def\E{{\math E}}
\def\P{{\math P}}
\def\Z{{\math Z}}
\def\Q{{\math Q}}
\def\C{{\math C}}
\def \e{{\rm e}}
\def \f{{\cal F}}
\def \g{{\cal G}}
\def \h{{\cal H}}
\def \d{{\tt d}}
\def \k{{\tt k}}
\def \i{{\tt i}}
\def \p{{\cal P}}
\def \s{{\cal S}}
\newcommand{\ed}{\mbox{$ \ \stackrel{d}{=}$ }}
\newtheorem{theorem}{Theorem}
\newtheorem{proposition}[theorem]{Proposition}
\newtheorem{lemma}[theorem]{Lemma}
\newtheorem{assum}[theorem]{Assumption}
\newtheorem{corollary}[theorem]{Corollary}
\centerline{\LARGE \bf Discretization methods}
\vskip 2mm
\centerline{\LARGE \bf  for homogeneous
fragmentations}

\vskip 1cm
\centerline{\Large \bf Jean Bertoin$^{(1)}$ and Alain Rouault$^{(2)}$}
\vskip 1cm
\noindent
\centerline{\bf September 2004}
\vskip 0.5cm
\noindent
(1) {\sl Laboratoire de Probabilit\'es et Mod\` eles Al\'eatoires
and Institut universitaire de France,
 Universit\'e Pierre et Marie Curie,  175, rue du Chevaleret,
 F-75013 Paris, France.}
\vskip 2mm
\noindent
(2) {\sl {LAMA, B\^atiment Fermat,
Universit\'e de Versailles, 45, avenue des Etats-Unis, F-78035
Versailles Cedex, France.
}}
\vskip 15mm

\noindent{\bf Summary. }{\small Homogeneous fragmentations describe the
evolution of a unit mass that breaks down randomly into pieces as time
passes. They can be thought of as continuous time analogs of a certain
type of branching random walks, which suggests the use of
time-discretization to shift known results from the theory of branching
random walks to the fragmentation setting. In particular, this yields 
interesting information about the asymptotic behaviour
of fragmentations.

On the other hand, homogeneous fragmentations can also be investigated
using a powerful technique of discretization of space due to
Kingman, namely, the theory of exchangeable partitions of
$\N$. Spatial discretization is especially well-suited to develop
directly for continuous times the conceptual method of probability
tilting of Lyons, Pemantle and Peres. }
 \vskip 5mm
\noindent
{\bf A.M.S. Classification.}  {\tt 60 J 25, 60 G 09. }
\vskip 3mm
\noindent{\bf e-mail.} {\tt $(1):$ jbe@ccr.jussieu.fr , $(2):$
rouault@math.uvsq.fr }

\section{Introduction}

Homogeneous fragmentations form a family of random processes in
continuous
times which have been introduced in
\cite{Be1}. Roughly, these are particle systems that model a
mass that breaks down randomly into pieces as time passes. More
precisely, each particle is identified with its mass (i.e.  it is
specified by a positive real number), and the fragmentation property
requires that different particles have independent evolutions.
The homogeneity property means that the process started from a single
particle with mass $x>0$ has the same
distribution as $x$ times the process started from a single particle
with unit mass.

This verbal description has obvious similarities with that of branching
random walks. More precisely, let us write $Z^{(t)}$ for the random
point
measure
which assigns a Dirac point mass at $\log x$ for every $x$ varying over
the set
of particles at time $t$. Taking logarithms transforms the fragmentation
and
homogeneity properties into the branching property for random point
measures.
More precisely, for every $t,t'\geq 0$, $Z^{(t+t')}$ is obtained from
$Z^{(t)}$ by replacing each atom $z=\log x$ of $Z^{(t)}$ by
a family $\{z+y, y\in{\cal Y}\}$, where ${\cal Y}$
is distributed as the family of the atoms of $Z^{(t')}$ for
$Z^{(0)}=\delta_0$,
and distinct atoms $z$ of $Z^{(t)}$ correspond to independent copies of
${\cal
Y}$.

Homogeneous fragmentations may be seen as extensions of so-called
branching random walks in continuous time. The latter have been
considered by Uchiyama \cite{Uch}, Biggins \cite{Bi3}, Kyprianou
\cite{Kyp1}, ... Their main feature
is that each particle has an exponentially distributed lifetime 
and at the instant of its death, scatters a random number of
children-particles in space relative to its death point
according to the point process. However the theory of branching processes
in continuous time does not encompass homogeneous fragmentations,
because usually each fragment starts to split instantaneously,
which would correspond to particles with zero lifetime in the branching
setting.

The close connection between homogeneous fragmentations and branching
random
walks suggests that one should try to reduce the study of fragmentations
to that of branching random walks by time-discretization. This is the path that we will
follow in the first part of this work.  Aside from some difficult technical problems (for
instance, a most useful notion such as the first branching time has no
analog for fragmentations since, in general, dislocations occur
instantaneously), this
enables us to shift several deep results on branching random walks to the
fragmentation setting. In particular, this yields 
interesting information about the asymptotic behaviour
of fragmentations, which refine earlier results in \cite{Be2}.

There is another discretization method that will play an important role in
this paper. The fundamental idea is due to
Kingman \cite{Ki}, who pointed out that partitions of an object, say with a unit
mass, can be fruitfully encoded by partitions of $\N$. In order to explain the
coding, we introduce a sequence of i.i.d. random points
$U_1,\ldots$ which are picked according to the mass distribution of the
object. One then considers at each time $t\geq0$ the random partition
of the set of indices
$\N$, such that two indices, say $i$ and
$j$, belong to the same block of the partition if and only if the points
$U_i$ and $U_j$ belong to the same fragment
of the object. By the law of large numbers, we see
that the masses of the fragments can be recovered as
the asymptotic frequencies of the blocks of the partition.
Roughly, the fundamental point in Kingman's coding is that it translates
the process to a 
discrete state-space. We refer to Pitman \cite{Pit} for
an important application of these ideas to a coalescent setting.

In the second part of this work, we shall present further applications of
Kingman's idea to homogeneous fragmentations. In particular, we will show
that spatial discretization is especially well-suited to adapt the
conceptual method of probability tilting introduced by Lyons, Pemantle and
Peres (see e.g. \cite{Lyo}) to homogeneous fragmentations,
which yields also some interesting limit theorems.

\section{Time discretization for ranked fragmentations}
This section is devoted to the presentation of some applications of 
 time-discretization to the asymptotic behaviour of ranked fragmentations
as time tends to $\infty$.
We shall first introduce some notation and definition; then we shall
merely translate results of Biggins on branching random walks in
continuous time in the special case when the so-called dislocation
measure of the fragmentation is finite. Finally, we shall extend the
preceding results to the case when the dislocation measure is infinite,
and also  derive from the extension of a result of Rouault sharp estimates for the probability of presence of
abnormally large fragments.

\subsection{Some notation and definition}
Throughout this paper, we will work
with the space of numerical sequences
$$\s\,:=\,\left\{{\bf s}=(s_1,\ldots): s_1\geq s_2\geq\cdots \geq0\hbox{ and }
\sum_{1}^{\infty}s_i\leq1\right\}$$
endowed with the uniform distance, which is a compact set.
A configuration $s\in\s$ should be thought of as the ranked masses of the
fragments resulting from the split of some object with unit total mass. 

We consider a family of Feller processes $X=(X_t, t\geq0)$ with
values in $\s$ and c\` adl\` ag paths.
For every $a\in[0,1]$, we let $\P_a$ denote the law
of $X$ with initial distribution $(a,0,\ldots)$ (i.e. the process
starts from a single fragment with mass $a$). We say that $X$ is a
 (ranked) {\it homogeneous fragmentation} if the following two
properties hold: 

\noindent $\bullet$ (Homogeneity property) For every $a\in[0,1]$, the law
of 
$aX$ under $\P_1$ is $\P_a$.

\noindent $\bullet$ (Fragmentation property) For every ${\bf s}=(s_1,
\ldots)\in\s$, the process started from $X(0)={\bf s}$ can be obtained as
follows. Consider $X^{(1)},
\ldots$ a sequence of independent processes with respective laws
$\P_{s_1}, \ldots$, and for every $t\geq0$, let $\hat X(t)$ be the random
sequence obtained by ranking in decreasing order the terms of the
random sequences $X^{(1)}(t), \ldots$. Then $\hat X$ has the law of $X$
started from ${\bf s}$.
 
It has been shown in \cite{Ber} and \cite{Be1} that homogeneous
fragmentations result from the combination of two different phenomena: a
continuous erosion and sudden dislocations. The erosion is a continuous
deterministic mechanism; dealing with
erosion is straightforward, and therefore we
will only consider homogeneous fragmentations with no erosion in this
work.

The dislocations occur randomly and can
be viewed as the jump-component of the process.
Roughly speaking, their distribution can be characterized by a 
measure $\nu$ on $\s$, called the {\it dislocation measure}. Informally
$\nu$ specifies the rates at which a unit mass splits; see the
forthcoming Section 2.2. It has to fulfil the conditions
$\nu(\{(1,0,\ldots)\})=0$ and 
\begin{equation}\label{clev}
\int_{\s}\left(1-s_1\right)\nu(d{\bf s})\,<\, \infty\,.
\end{equation}
More precisely, (\ref{clev}) is the necessary
and sufficient condition for a measure $\nu$ on $\s$ to be the dislocation
measure of some homogeneous fragmentation (see \cite{Ber} and \cite{Be1}).
We shall assume throughout this work that
\begin{equation}\label{pasperte}
\nu\left(\left\{{\bf s}\in\s:
\sum_{i=1}^{\infty}s_i<1\right\}\right)\,=\,0\,,
\end{equation}
 which means that no
mass is lost when a sudden dislocation occurs, and more precisely,
entails that the total mass is a conserved quantity for the fragmentation
process (i.e. $\sum_{i=1}^{\infty}X_i(t)=1$ for all $t\geq0$, 
$\P_1$-a.s.). 
In the sequel, we shall also implicitly exclude the
trivial case when $\nu\equiv 0$.  

Given a real number $r>0$, we say that a dislocation measure $\nu$ is
$r$-geometric  if $\nu$ is finite and is carried by the 
subspace of configurations ${\bf s}=(s_1,\ldots)\in\s$ such that
$s_i\in \{r^{-n}, n\in\N\}$. This holds if and only if
$\P_1(X_i(t)\in\{r^{-n}, n\in\N\}
\hbox{ for every }i\in\N)=1$ for all $t\geq0$.
We say that a dislocation measure 
is {\it non-geometric} if it is not $r$-geometric for any $r>0$.

We now introduce analytic quantities defined in terms of $\nu$
which will have an important role in this work.
First, we set
$$\underline p\,:=\,\inf\left\{p\in\R:
\int_{\s}\sum_{i=2}^{\infty}s_i^{p+1}\nu(d{\bf s})<\infty\right\}
\,,$$ and then for every $q>\underline p$
\begin{equation}\label{defphi}
\Phi(q)\,=\,\int_{\s}
\left(1-\sum_{i=1}^{\infty}s_i^{q+1}\right)\nu(d{\bf s})\,.
\end{equation}
The function $\Phi$ is a concave analytic increasing function; it is easy
to see (cf. Lemma 1 in \cite{Be2}) that the equation
$$\Phi(q)\,=\,(q+1)\Phi'(q)\,,\qquad q>\underline p$$
has a unique solution, which we denote by $\bar p$. More precisely, we
then have
\begin{equation}\label{defbarp}
\Phi(q)-(q+1)\Phi'(q) \,<\, 0 \quad
\Longleftrightarrow\quad  q\in]\underline p,
\bar p[\,,
\end{equation}
and 
\begin{equation}\label{eqcroiss}
\hbox{the map $q\to\Phi(q)/(q+1)$ increases on
$]\underline p,
\bar p[$ and decreases on $]\bar p,\infty[$.}
\end{equation}

To start with, we consider the simple sub-family of fragmentation
processes with a finite dislocation measure. These can be reduced to
continuous time branching random walks, and we specify some important
results of Biggins \cite{Bi3} in this setting.
Then we shall investigate the case when the dislocation measure is
infinite by time discretization.

\subsection{The case when the dislocation measure is finite}

Throughout this section, we assume that the dislocation measure $\nu$ is
finite. It is easy to construct a fragmentation process $X=(X(t),
t\geq0)$ with dislocation measure $\nu$.
Let the process start, say from the state ${\bf 1}:=(1,0,\ldots)$, and
stay there  for an exponential time with parameter $\nu(\s)$. Then the
process jumps independently of the waiting time to some random state in
$\s$ distributed according to the probability measure
$\nu(\cdot)/\nu(\s)$. After this first split, each fragment has a
similar evolution, independently of the 
other fragments. In words, a
fragment with mass
$x\in]0,1[$ breaks after some exponential time with parameter $\nu(\s)$,
and produces a random sequence of smaller fragments, say $xS$, where $S$
is a random variable in 
$\s$ with law $\nu(\cdot)/\nu(\s)$.

Plainly, the
empirical measure of the logarithms of the  fragments
\begin{equation}\label{defZ}
Z^{(t)}\,:=\,\sum_{i=1}^{\infty}\delta_{\log
X_i(t)}\,,\qquad t\geq0
\end{equation}
can be viewed as a branching random walk in continuous-time; see
 Uchiyama \cite{Uch}, Biggins \cite{Bi3}, Kyprianou
\cite{Kyp1}, ...
In this direction, let us identify two key quantities related to
branching random
walks in the fragmentation setting. 

First, it is easy
to see by an application of the
Markov property at the first splitting (i.e. branching) time, that the
Laplace transform of the intensity of the point process
$Z^{(t)}$ is given for $\theta>\underline p+1$ by
\begin{equation}\label{defm}
m(\theta)^t\,:=\,
\E\left(\int_{\R}\e^{\theta x}Z^{(t)}(dx)\right)\,
=\,\E\left(\sum_{i=1}^{\infty}X_i(t)^{\theta}\right)\,=\,
\exp(-t\Phi(\theta-1))\,.
\end{equation}
Second, there is also the identification for the so-called additive
martingale
$$
W^{(t)}(\theta)\,:=\,m(\theta)^{-t}\int_{\R} \e^{\theta
x}Z^{(t)}(dx)\,
=\,\exp(t\Phi(\theta-1))\sum_{i=1}^{\infty}X_i(t)^{\theta}
\,.$$

These observations allow us to apply Theorem 6 of Biggins
\cite{Bi3} (see also \cite{Bi1}), and we now state:

\begin{proposition}\label{P1}
Assume that the dislocation measure $\nu$ is finite. Then for every
$p>\underline p$, the process
$$M(p,t)\,:=\,W^{(t)}(p+1)\,=\,
\exp(t\Phi(p))\sum_{i=1}^{\infty}X_i^{p+1}(t)\,,\qquad
t\geq0$$ is a martingale with c\` adl\` ag paths. This martingale
converges uniformly on any compact subset of $]\underline p, \bar p[$, 
almost surely and in mean, as
$t\to\infty$. 
\end{proposition}
We also point out that for 
$p\geq\bar p$, the result of Biggins \cite{Bi1} entails that 
$\lim_{n\to\infty}M(p,n)=0$ a.s., and thus by the convergence theorem
of c\` adl\` ag nonnegative martingales, it holds that
$\lim_{t\to\infty}M(p,t)=0$ a.s.

\proof Fix some
compact interval
$[a,b]\subset ]\underline p,
\bar p[$, and recall that the event that
\begin{equation}\label{pvarfi}
\sum_{i=1}^{\infty}X_i^{a+1}(t)\,<\,\infty\qquad
\hbox{for all }t\geq0
\end{equation}
has probability one. In the sequel, we shall always work on this event.
We consider  the random continuous function on $[a,b]$
$$M(t): p\to \exp(t\Phi(p))\sum_{i=1}^{\infty}X_i^{p+1}(t)\,,$$
which defines a martingale with values in the Banach space ${\cal
C}([a,b],
\R)$.

Next, since $\sum_{i=1}^{\infty}s_i= 1$ for $\nu$-a.e. ${\bf s}\in\s$, observe
that, for every $\theta>\underline p+1$, we have
$\gamma(\theta-1)>\underline p$ for some $\gamma>1$ and then
Jensen's inequality implies that 
$$\int_{\s}\left(\sum_{i=1}^{\infty}s_i^\theta\right)^\gamma\nu(d{\bf s})
\,\leq\,\int_{\s}\left(\sum_{i=1}^{\infty}s_i^{\theta\gamma+1-\gamma}
\right)\nu(d{\bf s})\,<\,\infty\,.$$
On the other hand, thanks to (\ref{eqcroiss}), we get that
for every $\theta\in]\underline p+1, \bar p+1[$, there exists $\alpha\in
]1,\gamma[$  such that
$$m(\alpha\theta)<m(\theta)^\alpha\,.$$
By the argument used in the proof of Theorem 6 in Biggins \cite{Bi3}
(see also Remark (1) at the end of the present proof), all that needs to
be checked is the martingale $(M(t), t\geq0)$ has right-continuous paths.
We stress that the argument we give does not rely on the assumption
of finiteness for the dislocation measure. 

In this direction, it is convenient to use an
interval representation of the fragmentation (see \cite{Be3}).
Specifically, we consider the space of open subsets of the unit
interval, endowed with the metric induced by the Hausdorff distance for
the complementary closed set. One can construct a right-continuous family
$(\Theta(t), t\geq0)$ of random open subsets of the unit interval
such that for each $t\geq0$, $X(t)$ is the ranked sequence of the lengths
of the interval components of $\Theta(t)$, and $\Theta(t')\subseteq
\Theta(t)$ whenever $t\leq t'$. Recall also that the assumption
(\ref{pasperte}) ensures that each $\Theta(t)$ has full Lebesgue measure
a.s. 

For every  $x\in ]0,1[$ and $q\in[a,b]$, write 
$f_{t,q}(x)=|I_x(t)|^{q}$, where $I_x(t)$ denotes the interval
component of $\Theta(t)$ that contains $x$ and $|I_x(t)|$ its length.
In this setting, we thus have
$$\sum_{i=1}^{\infty}X_i^{q+1}(t)\,=\,\int_{0}^{1}f_{t,q}(x) dx\,.$$
For every $t_0\geq0$ and $x\in\Theta(t_0)$, we get from the
right-continuity of
$(\Theta(t), t\geq0)$ that  $|I_x(t)|$
increases to $|I_x(t_0)|$ as $t$ decreases to $t_0$.
The obvious
 upper-bounds
$$f_{t,q}(x)\,\leq\, f_{t',a}(x)\,,\qquad t\leq t'\hbox{ and }x\in]0,1[ \ \ q \leq 0$$
$$f_{t,q}(x)\,\leq\, f_{t_0,a}(x)\,,\qquad t\geq t_0\hbox{ and }x\in]0,1[\ \ q \geq 0$$
combined with (\ref{pvarfi})
enable
 us to apply the theorem of dominated convergence, and hence
$$\lim_{t\to t_0+}\int_{0}^{1}f_{t,q}(x)
dx\,=\,\int_{0}^{1}f_{t_0,q}(x) dx\,.
$$
This shows that with probability one, the real-valued martingales
$M(q,\cdot)$ have right-continuous paths for all $q\in[a,b]$.

To conclude, we observe that the random function
$q\to\sum_{i=1}^{\infty}X_i^{q+1}(t)$ is continuous and decreases as $q$
increases. An appeal to Dini's Theorem now shows that $(M(t), t\geq0)$
is right-continuous at all $t\geq0$, with probability one.
\QED

\noindent {\bf Remarks.} (1) It seems that the discretization argument in
the proof of the almost sure convergence in Theorem 6 in
\cite{Bi3} might require a further explanation. Indeed, it is observed
there that $W^{(n\delta)}$ converges a.s. when $n\to\infty$ through
integers for any $\delta>0$, and then claimed that this implies the
a.s. convergence of $W^{(t)}$ when $t\to\infty$ through the rationals.
The latter assertion does not look obvious, so we propose a slightly
different argument. One works with the martingale $W^{(t)}$ with values in
the space of continuous functions on some compact space, endowed with the
supremum norm, $||\cdot||$. We know that this martingale converges in
mean as $t\to\infty$. 
The norm is a convex map, therefore
for every integer $n$, the process
$||W^{(t+n)}-W^{(n)}||$ is a nonnegative submartingale with regular
paths. Doob's maximal inequality now entails that for every
$\varepsilon>0$
$$\P\left(\exists t>0: ||W^{(t+n)}-W^{(n)}||>\varepsilon\right)
\,\leq
\varepsilon^{-1}\sup_{t\geq0}\E\left(||W^{(t+n)}-W^{(n)}||\right)\,,$$
and the right-hand side converges to $0$ as $n\to\infty$.
An application of the Borel-Cantelli lemma now completes the proof of the
almost sure convergence of
$W^{(t)}$ as $t\to\infty$.

\noindent (2) We also mention that Proposition \ref{P1} can be extended
to complex numbers $p$ with $\underline p < \Re p < \bar p$.
The key point is to check 
that the martingale $p\to M(p,t)$ ($t\geq0$), viewed as a process with values in
the space ${\cal C}(K,\C)$ for some compact set $K\subset \{z: \underline p <
\Re p < \bar p\}$, has right-continuous paths. 
We know from the proof of Proposition \ref{P1} that the latter holds 
when $K$ is a real segment, and the general case follows easily.

\smallskip

Next, we derive an important
consequence of the preceding analysis concerning almost sure large
deviations for the empirical measure; see Corollary 4 and the discussion
on page 150 in  Biggins \cite{Bi3}.

\begin{corollary}\label{c1}
 Assume that 
the dislocation measure $\nu$ is finite and non-geometric. For 
 $p\in]\underline p, \bar p[$, let $M(p,\infty)$ be the terminal
value of the uniformly integrable martingale $M(p,\cdot)$.
 If $f:\R\to\R$ is a function with compact
support which is directly Riemann integrable,
then
$$\lim_{t\to\infty}  \sqrt t \, 
\e^{-t((p+1)\Phi'(p)-\Phi(p))}\int_{\R}^{}f(t
\Phi'(p)+y)Z^{(t)}(dy)\,=\, {M(p,\infty)\over \sqrt{2\pi
|\Phi''(p)|}}\int_{-\infty}^{\infty}f(y)\e^{-(p+1)y}dy\,.
$$
uniformly for $p$ in compact subsets of $]\underline p, \bar p[$, almost
surely.
\end{corollary}
We point out that this result applied for indicator functions of bounded
intervals gives a sharp large deviation statement that extends
Corollary 2 of \cite{Be2}.

\subsection{The case when the dislocation measure is infinite}
We now drop the assumption of finiteness of the dislocation measure
$\nu$, and merely assume that (\ref{clev}) holds.
As above we associate to the fragmentation $X$ the empirical
measures $Z^{(t)}$ defined in (\ref{defZ}); note that when
$\nu(\s)=\infty$, the process
$(Z^{(t)}, t\geq0)$ is no longer of the kind considered by Uchiyama \cite{Uch}.

\begin{theorem}\label{T1}
Proposition \ref{P1} and Corollary \ref{c1} hold  when one relaxes the
requirement of finiteness of the dislocation measure $\nu$  and merely
assumes (\ref{clev}).
\end{theorem}

\proof The
restriction of the process of empirical distribution to integers times
$(Z^{(n)}, n=0, 1,
\ldots)$ is a branching random walk; we aim at applying results of
Biggins \cite{Bi3} in this setting. 

The main difficulty is that the distribution of
$Z^{(1)}$ is not explicitly known in terms  of the dislocation measure
$\nu$. However, it is known that (\ref{defm}) still holds;
see e.g. the identity (6) in \cite{Be2}. Moreover, the proof of Theorem 2
 in \cite{Be2} shows that for every $p>\underline p$, there exists some
$\gamma>1$ such that for every $t\geq0$
\begin{equation}\label{momgfi}
\E\left(\left(\sum_{i=1}^{\infty}X_i^{p+1}(t)\right)^\gamma\right)
\,<\,\infty\,.
\end{equation}
Theorem 2 of Biggins \cite{Bi3} now shows that the conclusions
of Proposition \ref{P1} hold provided that $t\to\infty$ through integers.
We can complete the argument as in Theorem 6 of Biggins \cite{Bi3}
(using also the proof of Proposition \ref{P1} in the present
paper and the remark thereafter).
The extension of Corollary \ref{c1} is proven by adapting the arguments
of Biggins \cite{Bi3} on page 150. \QED

In a different direction, one can also use the skeleton method to
estimate the probability of presence of abnormally large fragments as
time goes to infinity. Indeed, a similar problem has been solved for branching
random walks, see \cite{Rou} and the references therein. Informally in
the so-called sub-critical region, the probability of presence of particles is
of the same order as the mean number of particles in that region, and the
asymptotic behaviour of the latter can be derived from the local central limit
theorem. This incites us to fix two real numbers
$\alpha<\beta$ and to introduce for every $t\geq0$ and $x\in\R$ the
notation:
\begin{eqnarray}
U(t, x) &:=& \P( Z^{(t)}
\big([x + \alpha, x + \beta]\big) > 0) 
\\
V(t, x) &:=& \E \left(Z^{(t)} 
\big([ x + \alpha  ,  x + \beta ]\big)\right)\,.
\end{eqnarray}

\begin{theorem} Assume that the dislocation measure $\nu$ is
non-geometric. 

\label{presfrag}
\begin{itemize}
\item[{\rm (i)}] If $p > \underbar p$, we have 
\begin{eqnarray*}
\lim_{t\to\infty} 
\sqrt t \, \e^{-t((p+1)\Phi'(p)-\Phi(p))} 
V(t, -t\Phi'(p))  \,=\, \frac{1}{\sqrt{2\pi
|\Phi''(p)|}} 
(p+1)^{-1}
 \left(\e^{-(p+1) \alpha}-\e^{-(p+1) \beta}\right)
\,. 
\end{eqnarray*}
\item[{\rm (ii)}] If $p > \bar p$, there exists a positive finite constant
$K_{p}$ such that
\begin{eqnarray*}
\label{usurvfg}
\lim_{t \rightarrow \infty}\frac{U(t, -t\Phi'(p))}{V(t, -t\Phi'(p) )} = K_p \,. 
\end{eqnarray*}
\end{itemize}
\end{theorem}
We point out that in the range $p \in ]\underbar {\it p} , \bar p[$,  $(i)$ is the
counterpart in mean of the result of Corollary 2, when $f = {\bf 1}_{[\alpha , \beta]}$. 
It implies that the convergence there holds also in $L^1(\P)$, thanks to Scheff\'e's
theorem. 
\smallskip

\proof The proof relies on a result of Rouault \cite{Rou} 
for branching random walks, which has been extended under more general
conditions in \cite{pageweb}, and that we take here for granted.
Recall also (\ref{momgfi}), which ensures that 
the assumptions of \cite{pageweb} are satisfied.

We apply the
skeleton method. Let
$h > 0$ be a time mesh; the fragmentation process observed at times 
$nh$ ($n \in \N$) yields a branching random walk. Write
$$\widehat{Z^{(h)}}(\theta):=\int_{\R} \e^{\theta x} Z^{(h)}(dx)=
\sum_i  X_i (h)^\theta
$$ 
and 
$$\Lambda_h(\theta) := \log \E\left(\widehat{Z^{(h)}}(\theta)\right)
= - h\Phi(\theta-1)\,.$$

In the case $h=1$, we write for simplicity $\Lambda=\Lambda_1$
and set
$$a =  - \Phi'(p)\ ,\ \sigma_p ^2 = -\Phi'' (p)\ ,\ 
\Lambda^*(a)=\Phi(p)-(p+1)\Phi'(p)\,.$$
For an arbitrary mesh $h>0$, we define by scaling
$$\Lambda_h^* (x) = h \Lambda^* (x/h)\,.$$
It is immediately checked that if  $\theta$ solves $a=\Lambda_h'(\theta)$,
then $\Lambda_h^* (a)=\theta \Lambda_h'(\theta)-\Lambda_h(\theta)$.
Applying Theorem 2 
 in \cite{pageweb}, we get first
\begin{eqnarray}
\label{ckv}
\lim_n \sigma_p \sqrt{2\pi n h}\ \e^{nh\Lambda^* (a)}V(nh, anh)  
= (p+1)^{-1}\left(\e^{-(p+1)\alpha} -\e^{-(p+1)\beta}\right),
\end{eqnarray}
 and then that, if $\Lambda^*_h (a) > 0$ the limit
$$\lim_n \frac{U(nh, anh)}{V(nh, anh)} 
=: K_p^{(h)} $$
exists for each $h>0$ and is positive.
Now, it is easy to see that  the functions
\begin{eqnarray}
t \mapsto \sigma_p \sqrt{2\pi t}\ \e^{ t\Lambda^* (a)} V(t, at)\ \
\hbox{and} \ \ t \mapsto \frac{U(t, at)}{V(t, at)} 
\end{eqnarray}
are continuous. We  apply the Croft-Kingman lemma (\cite{AsHer} A 9.1 p.438,
see also \cite{King}).  Both limits exist  when $t\rightarrow \infty$. In
the first case, it is of course the right hand side of (\ref{ckv}). In
the second case, it is any $K_p ^{(h)}$ since they are all equal.
\QED

\section{Spatial discretization and fragmentation of partitions}
This section is devoted to another useful discretization technique which has been sketched in the Introduction
and that we now recall for convenience.
We may suppose that we are given a
sequence of i.i.d. random points $U_1,\ldots$ which are picked according
to the mass distribution of the object. These random points are assumed
to be independent of the fragmentation process. One then looks at the
fragmentation process as the evolution as time $t$ passes of the random
partition $\Pi(t)$ of $\N$ which is given as follows. Two indices $i$
and $j$ are in the same block of $\Pi(t)$ if and only if the points
$U_i$ and $U_j$ belong to the same fragment
of the object at time $t$. Plainly, the random partition gets finer and finer
as time passes.

We first present in Section 3.1 the necessary background on partitions of $\N$,
and then in Section 3.2 the Poissonian construction of homogeneous
fragmentations and the connection with subordinators, following closely
\cite{Be1}.
Section 3.3 is devoted to the study of probability tilting based on
additive martingales, adapting to the random partition setting the
so-called conceptual method of Lyons et al. \cite{Lyo, LPP}. Finally, 
as an example of application, we investigate in Section 3.4 the
convergence of the so-called derivative martingale.

\subsection{Preliminaries}
A {\it partition} of $\N=\{1,\ldots\}$ is a sequence
$\pi=\left(\pi_1,\pi_2,\ldots\right)$ of disjoint subsets, called  {\it
blocks},
such that $\bigcup \pi_i=\N$.  The blocks $\pi_i$ of a partition are
enumerated
in the increasing order of their least element, i.e.
$\min \pi_i\leq \min \pi_j$ for $i\leq j$, with the convention that
$\min \emptyset =\infty$.
If $\pi$ and $\pi'$ are two partitions of $\N$, we say that $\pi$ is
finer
than $\pi'$ if every block of $\pi$ is contained into some block of
$\pi'$.

For every block $B\subseteq
\N$, we denote by
$\pi_{\mid B}$ the partition of $B$ induced by $\pi$ and an obvious
restriction. For every integer $k$, the block
$\{1,\ldots,k\}$ is denoted by $[k]$. A partition $\pi$
is entirely determined by the sequence of its restrictions
$\left(\pi_{\mid [k]},k\in\N\right)$, and conversely, if for every
$k\in\N$,
$\gamma_k$ is a partition of $[k]$ such that the restriction
of $\gamma_{k+1}$ to $[k]$ coincides with $\gamma_k$ (this will be
referred
to the {\it compatibility property} in the sequel), then there exists a
unique partition $\pi\in\p$ such that $\pi_{\mid [k]}=\gamma_k$ for every
$k\in\N$.

The space of partitions of $\N$ is denoted by $\p$ and
endowed with the hyper-distance
$${\rm dist}(\pi,\pi')\,:=\,1/\max\left\{k\in\N: \pi_{\mid[k]}=
\pi'_{\mid[k]}\right\}\,,$$
with the convention $1/\max \N := 0$.
This makes $\p$ compact.

One says that a block $B\subseteq \N$ has an {\it asymptotic frequency}
if the limit
$$|B|\,:=\,\lim_{n\to\infty}n^{-1}{\rm Card}(B\cap [n])$$
exists. When every block of some partition $\pi\in\p$ has an asymptotic
frequency, we write
$|\pi|=(|\pi_1|,\ldots)$, and then
$|\pi|^{\downarrow}=(|\pi|^{\downarrow}_1,\ldots)\in\s$ for the decreasing
rearrangement\footnote{Ranking the asymptotic frequencies of the blocks
of
$\pi$ in the decreasing order is just a simple procedure to forget the
labels of these blocks. In other words, we want to consider the family
of the asymptotic frequencies without keeping the
additional information provided by the way blocks are labelled.}
of the sequence $|\pi|$.
In the case when some block of the partition $\pi$ does not have an
asymptotic frequency, we decide that
$|\pi|=|\pi|^{\downarrow}=\partial$,
where $\partial$ stands for some extra point added to $\s$.
This defines a natural map  $\pi\to|\pi|^{\downarrow}$ from $\p$ to
$\s\cup \{\partial\}$ which is not continuous.

We call {\it nested partitions} a collection $\Pi=\left(\Pi(t),
t\geq0\right)$ of partitions of $\N$ such that
$\Pi(t)$ is finer than $\Pi(t')$ when $t'\leq t$.
There is a simple procedure for the construction of a large family of
nested partitions which we now describe and will use throughout the rest
of this section.

 We call {\it discrete point measure} on
$\R_+\times
\p\times \N$ any measure $\omega$ which can be expressed in the form
$$\omega\,=\,\sum_{(t,\pi,k)\in {\cal D}}^{\infty}\delta_{(t,\pi,k)}$$
where ${\cal D}$ is a subset of $ \R_+\times \p\times \N$ such that
the following two requirements hold:

\noindent $\bullet$ For every $t\in\R$, $\omega(\{t\}\times \p\times \N)=0$ or $1$.

\noindent $\bullet$  For
every real number $t'\geq0$ and integer
$n\geq1$
$${\rm Card}\left\{(t,\pi,k)\in {\cal D}: t\leq t', {\pi}_{\mid
[n]}\neq {\rm trivial}(n),k\leq n\right\}\,<\,\infty\,,$$
where ${\rm trivial}(n)=\left([n],\emptyset, \emptyset,\ldots\right)$
stands
for the partition of $[n]$ which has a single non empty block\footnote{Roughly, ${\rm
trivial}(n)$ plays the role of a neutral element 
in the space of partitions of $[n]$.}.

Starting from an arbitrary discrete point measure $\omega$ on
$\R_+\times \p\times \N$, we may construct
nested partitions $\Pi=(\Pi(t), t\geq0)$ as follows: Fix $n\in\N$; the
assumption that the point measure $\omega$ is discrete enables us to
construct a step-path $(\Pi(t,n), t\geq0)$ with values in the space of
partitions of $[n]$, which only jumps at times
$t$ at which the fiber $\{t\}\times\p\times\N$ carries an atom of
$\omega$,
say $(t,\pi,k)$,
such that $\pi_{\mid [n]}\neq {\rm trivial}(n)$ and $k\leq n$. In that
case,
$\Pi(t,n)$ is the partition obtained by replacing the $k$-th block of
$\Pi(t-,n)$, viz. ${\Pi_k(t-,n)}$, by the restriction $\pi_{\mid
{\Pi_{k}(t-,n)}}$ of $\pi$ to this block, and leaving the other blocks
unchanged.
Now it is immediate from this construction that for each time $t\geq0$,
the
sequence
$(\Pi(t,n), n\in
\N)$ is compatible, and hence there exists a unique partition $\Pi(t)$
such that $\Pi(t)_{\mid [n]}=\Pi(t,n)$ for each $n\in\N$.

\subsection{Poisson measures, homogeneous fragmentations, and subordinators}

We denote the space of discrete point measures on $\R_+\times
\p\times \N$ by $\Omega$, and the
sigma-field generated by the restriction to
$[0,t]\times
\p\times \N$ by
$\g(t)$. So $(\g(t))_{t\geq0}$ is a filtration, and the nested
partitions
$(\Pi(t), t\geq0)$ are $(\g(t))_{t\geq0}$-adapted. We shall also need to
consider
the sigma-field
$\f(t)$  generated by the decreasing
 rearrangement $|\Pi(r)|^{\downarrow}$
of the sequence of the asymptotic frequencies of the blocks of $\Pi(r)$
for
$r\leq t$, and $(\f(t))_{t\geq0}$ is a sub-filtration
of  $(\g(t))_{t\geq0}$.

Now consider a dislocation measure
$\nu$, i.e. a measure on $\s$ which fulfils the requirements of Section
2.1. According to Theorem 2 in \cite{Be1}, there exists a
unique measure
$\mu$ on $\p$ which is {\it exchangeable} (i.e. invariant by the action of
finite permutations on $\p$), and such that $\nu$ is
the image of $\mu$ by the map $\pi\to|\pi|^{\downarrow}$.
An important fact which stems from exchangeability, is that the
distribution of
the asymptotic frequency of the first block $|\pi_1|$ under the measure
$\mu$ is that of a size-biased picked term from the ranked sequence $s$
under
$\nu$. In other words, there is the identity
\begin{equation}\label{sizebmu}
\int_{\p}f(|\pi_1|)\mu(d\pi)\,=\,
\int_{\s}\sum_{i=1}^{\infty}s_if(s_i)\nu(d{\bf s})\,,
\end{equation}
where $f:[0,1]\to\R_+$ denotes a generic measurable function with
$f(0)=0$.

Let $\P$ be the probability
measure on $\Omega$ corresponding to the law  of a Poisson point measure
with intensity $dt\otimes \mu\otimes \#$,
where $\#$ denotes the counting measure on $\N$.
The assumption (\ref{clev}) on the dislocation measure $\nu$ ensures that
$\omega$ is a discrete point
measure $\P$-a.s.
The nested partitions $\left(\Pi(t), t\geq0\right)$ constructed above 
from $\omega$  now form a  Markov process; see Section 3 in \cite{Be1}.
More precisely the  Markov property is essentially a variation of the
branching property; it  can be stated as follows. For every
$t,t'\geq0$, the conditional distribution of $\Pi(t+t')$ given $\g(t)$
is
the
same as that of the random partition of $\N$ induced by the
restrictions $\Pi^{(1)}(t')_{\mid B_1},
\Pi^{(2)}(t')_{\mid B_2},\ldots$, where
$\Pi^{(1)}, \ldots$ are independent
copies of $\Pi$ and $(B_1,\ldots)=\Pi(t)$ is the sequence of
blocks of $\Pi(t)$. In the terminology of \cite{Be1}, we say that
$\Pi=\left(\Pi(t), t\geq0\right)$ is a (partition valued) {\it homogeneous
 fragmentation} under
$\P$. 

Another crucial fact is that the partitions $\left(\Pi(t), t\geq0\right)$
are exchangeable under $\P$, i.e. their distribution is 
invariant under the action of finite
permutations on $\N$; see Section 3 in \cite{Be1}.
It follows from a celebrated theorem of Kingman \cite{Ki} that $\P$-a.s., $\Pi(t)$ has
asymptotic frequencies for all
$t\geq0$; cf. Theorem 3(i) in \cite{Be1}.
The process of ranked asymptotic frequencies $|\Pi|^{\downarrow}:=X$
is a Markov process with values in $\s$; it provides a version of the
ranked fragmentation which we considered in Section 2; cf. \cite{Ber}.

The {\it tagged fragment} is the fragment of the object
that contains the first tagged point $U_1$, i.e. which corresponds to the
first block $\Pi_1(\cdot)$. 
The process $|\Pi_1(\cdot)|$ of the asymptotic frequencies of the first
block and
its logarithm,
$$\xi(t)\,:=\,-\log|\Pi_1(t)|\,,\qquad  t\geq0$$
 will have a special role
in this study. A crucial point is that
 under $\P$, $\xi=(\xi_t,
t\geq0)$ is a subordinator with Laplace exponent $\Phi$, which is given by
(\ref{defphi}); cf. Theorem 3(ii) in \cite{Be2}. 
This means that $(\xi(t), t\geq0)$ is a c\` adl\` ag process with
independent and
stationary increments, and the Laplace
transform
of its one-dimensional distribution is given by the identity
$$\E\left(\exp(-q \xi(t)\right)\,=\,\exp(-t\Phi(q))\,,\qquad
q>\underline
p\,.$$

More precisely,
let us denote by $\g_1(t)$ the sigma-field generated by
the restriction of the discrete point measure $\omega$ to the fiber
$[0,t]\times \p\times \{1\}$. So $(\g_1(t))_{t\geq0}$ is a
sub-filtration
of  $(\g(t))_{t\geq0}$, and the first block of $\Pi(t)$, $\Pi_1(t)$, and
a fortiori its asymptotic frequency  $\e^{-\xi_t}$, are
$\g_1(t)$-measurable. Let ${\cal D}_1\subseteq [0,\infty[$ be the random
set
of times $r\geq0$ for which the discrete point
measure has an atom on the fiber $\{r\}\times \p\times \{1\}$, and for
every $r\in{\cal D}_1$, denote the second component of this atom by
$\pi(r)$. The construction of the nested partitions from the discrete
point
measure yields  the identity
\begin{equation}\label{eqC}
\exp\left(-\xi_t\right)\,=\,|\Pi_1(t)|\,=\,\prod_{r\in
{\cal D}_1\cap [0,t]}|\pi_1(r)|\,,
\end{equation}
for all $t\geq0$, a.s. under $\P$; see e.g. the first remark at the end
of Section 5 in \cite{Be1}.
Observe that taking logarithm turns the identity (\ref{eqC})
into the L\'evy-It\^o decomposition for subordinators.

Finally, the conditional distribution of the size of the tagged fragment,
$|\Pi_1(t)|=\e^{-\xi(t)}$, given $\f(t)$ (the sigma-field generated by the
ranked asymptotic frequencies) is
that of a size-biased sample from the ranked sequence
$|\Pi(t)|^{\downarrow}$. In other words, we have
$$\E\left(f(\exp(-\xi(t))\right)\,=\,\E\left(\sum_{i=1}^{\infty}
|\Pi_i(t)|f(|\Pi_i(t)|)\right)
\,=\,\E\left(\sum_{j=1}^{\infty}
X_j(t)f(X_j(t))\right)$$
where
$f:[0,1]\to\R_+$ denotes a generic measurable function
with $f(0)=0$.
More generally, exchangeability ensures that for every $t\geq0$, the
sequence
$|\Pi(t)|$ of the asymptotic frequencies is a size-biased reordering of
the
ranked sequence $X(t)=|\Pi(t)|^{\downarrow}$.

\subsection{Additive martingales and tilted
probability measures}

There are two simple martingales connected to fragmentations for every
parameter
$p>\underline p$ :
First, a well-known fact for subordinators is that
$${\cal
E}(p,t)\,:=\,\exp(-p\xi(t)+t\Phi(p))\,=\,\e^{t\Phi(p)}|\Pi_1(t)|^p$$
is a positive  $(\P,\g(t))$-martingale.
Second, when we project ${\cal E}(p,t)$ on the sub-filtration $\f(t)$,
we recover the additive martingale
$$M(p,t)\,=\,\exp(t\Phi(p))\sum_{i=1}^{\infty}|\Pi_i(t)|^{p+1}
\,=\,\exp(t\Phi(p))\sum_{j=1}^{\infty}X_j^{p+1}\,.$$
We point out that, more precisely, $M(p,\cdot)$ is a
$(\P,\g(t))$-martingale  which is
adapted to the sub-filtration
$\f(t)$. 

Following the genuine method of Lyons, Pemantle and Peres (see e.g.
\cite{Lyo}), we introduce the
{\it tilted probability measure} $\P^{(p)}$ on the space of discrete
point
measures
$\Omega$ endowed with the filtration $(\g(t))_{t\geq0}$ by
\begin{equation}\label{eq1}
d\P^{(p)}_{\mid
\g(t)}\,=\,{\cal E}(p,t) \, d\P_{\mid
\g(t)}\,.
\end{equation}
Observe that projections on the sub-filtration $\f(t)$ give the identity

\begin{equation}\label{eq2}
d\P^{(p)}_{\mid
\f(t)}\,=\,M(p,t) \, d\P_{\mid
\f(t)}\,.
\end{equation}
The effect of the change of probability is easy to describe, both at
the level of the tagged fragment and that of the discrete point
measure.
\begin{proposition}\label{P3} 
{\rm (i)} 
Under
$\P^{(p)}$, the process $\xi_t=-\log|\Pi_1(t)|$ is a subordinator with
Laplace exponent
$$\Phi^{(p)}(q)\,:=\,\Phi(p+q)-\Phi(p)\,,\qquad q>\underline p-p\,.$$

\noindent{\rm (ii)} 
Under
$\P^{(p)}$, the discrete point measure $\omega$ is Poissonian. More
precisely:

\noindent $\bullet$  The restriction of $\omega$ to  $\R_+\times \p\times
\{2,3,\ldots\}$ has the same distribution as under $\P$ and
is independent  of the
restriction to the fiber  $\R_+\times \p\times
\{1\}$.

\noindent $\bullet$ In the notation of
Section 3.1, the family
$\{(r, \pi(r)), r\in{\cal D}_1\}$ is that of the atoms of a Poisson
random measure on $\R_+\times \p$ with intensity $dr\otimes \mu^{(p)}$,
where
$$\mu^{(p)}(d\pi)\,=\,|\pi_1|^{p}\mu(d\pi)\,.$$
\end{proposition}

\proof The first assertion stems from the classical Esscher transform;
see for instance Example 33.15 in Sato \cite{Sa}.
The description of the law of the discrete point measure under
$\P^{(p)}$ is easily seen from the formula
(\ref{eqC}), and classical properties of the exponential tilting for
Poisson random measures.
\QED

We stress that the tilting only
affects the distribution of the discrete point measure on the fiber
$\R_+\times \p\times \{1\}$.  
In this direction, it is interesting to use Proposition \ref{P3}(ii) and compare the
evolution
of
the ranked-fragmentation $X=|\Pi(\cdot)|^{\downarrow}$  under
$\P^{(p)}$ with the evolution under $\P$. For the sake of simplicity, 
we again suppose here that the dislocation measure 
$\nu$ is finite, so the evolution of the ranked fragmentation under $\P$
is described in section 2.2.

First, we observe that by absolute continuity, the random
partition  $\Pi(t)$ obtained by evaluating the nested partitions at
time $t$, possesses asymptotic frequencies a.s. under the tilted
probability
$\P^{(p)}$. 
The first block $\Pi_1(\cdot)$ has a special role in the
definition
of the tilted probability $\P^{(p)}$, and cannot
be recovered from the ranked sequence $|\Pi(\cdot)|^{\downarrow}$ alone.
Let us call {\it marked} the unique particle (i.e. asymptotic frequency)
at time $t$ corresponding to
$\Pi_1(t)$ and {\it unmarked} the other ones.
Observe that under $\P$, the tagged fragment coincides with the marked
particle, and that under
$\P^{(p)}$, the unmarked particles follow the same evolution as under
$\P$,
i.e. they split according to $\nu$, independently of the others, and
only
produce
unmarked particles. Under $\P^{(p)}$, the marked particle splits
independently of the other
particles, but with a different rate, namely
$$\nu^{(p)}(ds)\,:=\,\left(\sum_{i=1}^{\infty}s_i^{p+1}\right)\nu(d{\bf s})\,.$$
Indeed, we deduce from Proposition \ref{P3}(ii) that $\nu^{(p)}$ is the
image of the intensity measure
$\mu^{(p)}$ by the map $\pi\to |\pi|^{\downarrow}$, and since under $\mu$,
$|\pi_1|$
can be
viewed as a size-biased pick from the ranked sequence
$|\pi|^{\downarrow}$ (recall the identity (\ref{sizebmu})), this
yields the formula above.  The ``new'' marked particle is picked at
random
amongst the particles produced by the splitting of the ``old'' marked
particle
as follows: Let $x$ denote the mass of the old marked particle and $xs$
the
ranked sequence of the masses of the particles produced after the
splitting,
where
$s=(s_1,\ldots)\in\s$. Then the probability that the new marked particle
has
mass $xs_j$ equals
$s_j^{p+1}/\sum_{i=1}^{\infty}s_i^{p+1}$.
In short, the marked particle $\Pi_1(\cdot)$ can  be viewed
as a canonic analog in the fragmentation setting of the so-called
{\it spine} in the branching random walk framework.

 It may be interesting to point at the following
connection with the so-called {\it thinning}  of discrete point measures.
Recall that, given some metric space $A$, a discrete point measure
$\mu$ on $A$ with atoms $x_1,\ldots$ (i.e. $\mu=\sum_{i=1}^{\infty}\delta_{x_i}$),
and a measurable function $f:A\to[0,1]$, an $f$-thinning of
$\mu$ is the random discrete point measure $\mu^{(r)}$ obtained by
keeping each atom $x$ of
$\mu$ with probability $f(x)$, independently of the others. In other words, 
$$\mu^{(f)}\,=\,\sum_{j=1}^{\infty}{\bf 1}_{\{ \chi_j=1\}}\delta_{x_j}$$
where $\chi_1,\ldots$ is a sequence of independent Bernoulli variables with
$\P(\chi_j=1)=f(x_j)$.
Informally, the following corollary shows that dislocations
are less (respectively, more) frequent under $\P^{(p)}$ than under 
$\P$ for $p>0$ (respectively for $p<0$).

\begin{corollary}\label{C1} For every $q>0$, let $f_q$ be the map on $\R_+\times \p\times \N$
such that $f_q(t,\pi,k)=|\pi_1|^{q}$ if the partition $\pi$ possesses asymptotic frequencies
and $k=1$, and $f_q(t,\pi,k)=0$ otherwise.

\noindent {\rm (i)} For every $p>0$, the image of $\P$ by an $f_p$-thinning
of the discrete point measure $\omega$ is $\P^{(p)}$.

\noindent {\rm (ii)} For every $p\in]\underline p,0[$, the image of $\P^{(p)}$
by an
$f_{-p}$-thinning of the discrete point measure $\omega$ is $\P$.
\end{corollary}

\proof The first statement is immediate from Proposition \ref{P3}(ii) and properties of
thinning  (see e.g. Chapter 5 in Kingman \cite{KP}). Suppose now that $\underline p<p<0$, and
work under $\P$. Let $\omega'$ be a random Poisson point measure on the fiber $\R_+\times
\p\times
\{1\}$ with intensity $dt\otimes (|\pi_1|^{p}-1)\mu(d\pi)$, which is independent of
$\omega$. By superposition of independent Poisson measures, $\P^{(p)}$
can be identified as the law of
$\omega+\omega'$ under $\P$. It follows readily that the original
probability measure $\P$ can be recovered from $\P^{(p)}$ by $f_{-p}$-thinning.
\QED

\subsection{The derivative martingale}

We end this work by considering the so-called
 derivative martingale that we now introduce. Recall that $\bar p>0$ is
the
critical value for the convergence in $L^1(\P)$ of the additive
martingales.
The
process
$${\cal E}'(t)\,:=\,\left(t\Phi'(\bar p)-\xi(t)\right)\exp(-\bar
p \xi(t)+t\Phi(\bar p))\,,\qquad t\geq0$$
is clearly a $(\P,\g(t))$-martingale; its projection on the
sub-filtration
$\left(\f(t)\right)_{t\geq0}$ is a $(\P,\f(t))$ martingale,
called the {\it derivative martingale} and given by
$$M'(t)\,=\,\sum_{i=1}^{\infty}
\left(t\Phi'(\bar p)+\log\left(|\Pi_i(t)|\right)\right)
\exp\left(t\Phi(\bar p)\right)|\Pi_i(t)|^{\bar p +1}\,.$$
We stress that the derivative martingale is not always positive, which
contrasts
with the case of additive martingales.
The idea of considering the derivative martingale at the critical value
goes back to Neveu \cite{Nev}
for the branching Brownian motion. For the branching random walk,
it has been considered by Kyprianou \cite{Kyp2}, Liu \cite{Liu} with
the help of a functional equation and by Biggins and Kyprianou
\cite{Bi4} with the measure change method. 

\begin{proposition}\label{P5}
 {\rm (i)}
The martingale
$M'$ converges $\P$-a.s. to a finite non-positive limit $M' (\infty)$,

\noindent {\rm (ii)}
$\E( M' (\infty)) = -\infty$,

\noindent {\rm (iii)}
$\P
(M' (\infty) < 0) = 1$.

\end{proposition}

The proposition could be derived
by time discretization from its analog for branching random walks.
However, as checking the technical details may be rather involved in this
instance, we shall present a direct proof based on probability tilting.
This technique is due to 
Lyons et al. \cite{LPP}; it can be also applied
to establish the uniform integrability of additive martingales (cf.
Theorem \ref{T1}). We also refer to Harris \cite{Ha} and Kyprianou
\cite{Kyp3} for related treatments.

{\noindent{\bf Proof of (i):}\hskip10pt}  Define
for every $i\in\N$ and $s \leq t$,
$\beta_{s,t} (i)$ as the unique block of $\Pi(s)$ containing $\Pi_i(t)$.
For $a>0$, let
$$
\cases { \Pi_i ^{(a)} (t) =\Pi_i (t), & if $|\beta_{s,t} (i)|\leq
\exp\{a-s\Phi'(\bar
p)\}$ for every $s\leq t$; \cr
\Pi_i ^{(a)} (t)= \emptyset, &otherwise . \cr}
$$
The family
$\left\{\Pi_i^{(a)} (t):  i\in\N \right\}$
obviously possesses asymptotic frequencies.
Moreover, it should be plain that as $t$ varies in $[0,\infty[$, this
family of partitions is nested. We denote by $\left({\cal
H}(t)\right)_{t\geq0}$ the
filtration generated by the process of their ranked asymptotic
frequencies,
so $\left({\cal H}(t)\right)_{t\geq0}$ is another sub-filtration of
$\left(\g(t)\right)_{t\geq0}$.

Because $\Phi'(\bar p)=\Phi(\bar p)/(\bar p +1)$ and the martingale
$M(\bar p,t)$ converges to $0$,
$\P$-a.s., we have   $\sup_{t\geq0}\left\{\exp(t\Phi'(\bar
p))|\Pi(t)|_1^{\downarrow}\right\}<\infty$, $\P$-a.s.
It follows that
$$\lim_{a\to\infty}\P\left(\Pi_i^{(a)} (t) = \Pi_i (t) \hbox{ for all
}i\in\N \ \ \hbox{ and for all} \ t \geq 0\right)
\,=\,1\,.$$
Thus, in order 
to prove the
existence of a finite limit for $M'$, it suffices to
establish that if
\begin{equation}
\label{defMa}
M_a (t) := \sum_{i=1}^{\infty}
\left(\log\left(1/|\Pi_i(t)|\right)-t\Phi'(\bar p)+a\right)
\exp\left(t\Phi(\bar p)\right)|\Pi_i^{(a)}(t)|^{\bar p+1}
\end{equation}
then $\lim_{t\rightarrow \infty} M_a (t) =: M_a (\infty)$ exists $\P$-a.s. for every $a>0$.

From now on, we fix $a>0$. Since the process $\xi(t)-\Phi'(\bar p)t$
has no negative jumps,
$${\cal M}_a(t) \,:=\,
\left(\xi(t)+a-t\Phi'(\bar p)\right)\exp(-\bar p\xi(t)+t\Phi(\bar p))
{\bf 1}_{\{t<\zeta_a\}}$$
where $\zeta_a = \inf\{t\geq0: \xi(t) <  t\Phi'(\bar p) -a \}$,
can be viewed as a stopped (non-negative) $\P$-martingale.
Its projection on the
sub-filtration
$\left({\cal H}(t)\right)_{t\geq0}$ is
$M_a (t)$, which therefore is
a non-negative $(\P,{\cal H}(t))$ martingale, and thus possesses a finite 
limit as $t \rightarrow \infty$, $\P$-a.s.  \QED

{\noindent{\bf Proof of (ii):}\hskip10pt}It is sufficient to show that for all $a>0$, in the
notation above,
\begin{equation}\label{eqliminf}
\liminf_{t\to\infty}M_a(t)<\infty\,, \qquad
\Q \hbox{-a.s.}
\end{equation} Indeed (\ref{eqliminf}) entails that the
$\P$-martingale $M_a$ is  uniformly  integrable (see Lyons
\cite{Lyo}).  
Then $\E \left(M_a (t)\right) = a$, which 
yields $\E( -M'(\infty) )\geq a$ for every
$a > 0$.
In this direction, we
introduce the tilted probability measure
$\Q$ on
$\Omega$ given by
\begin{equation}\label{eq5}
d\Q_{\mid
\g(t)}\,=\,a^{-1}{\cal M}_a(t) \, d\P_{\mid
\g(t)}\,;
\end{equation}
so we also have
$$d\Q_{\mid
{\cal H}(t)}\,=\,a^{-1} M_a(t) \, d\P_{\mid
{\cal H}(t)}\,.$$
The following lemma, which is a simple variation of Proposition
\ref{P3}, lies in the heart of the proof of (\ref{eqliminf}).

\begin{lemma}\label{P3'} 
{\rm (i)} 
Under
$\Q$, the process
$$\lambda(t)\,:=\,\left(\xi(t)+a-t\Phi'(\bar p)\right){\bf
1}_{\{t<\zeta_1\}}\,,\qquad t\geq0$$ is  a centered L\'evy process with no negative jumps started from $a$ and conditioned to stay
positive forever (see for instance \cite{Chau}).

\noindent{\rm (ii)} 
Under
$\Q$, the restriction of $\omega$ to  $\R_+\times \p\times
\{2,3,\ldots\}$ has the same distribution as under $\P$ and
is independent  of the
restriction to the fiber  $\R_+\times \p\times
\{1\}$.
\end{lemma}

We are now able to complete the proof of Proposition \ref{P5}(ii). 

{\noindent{\bf Proof of {\rm (\ref{eqliminf})}:}\hskip10pt}
It is easily seen from Lemma \ref{P3'}(i) that
\begin{equation}\label{eq6}
\inf\left\{\lambda(t),
t\geq0\right\}\,>\,0\quad \hbox{and}\quad
\lim_{t\to\infty}{\log \lambda(t)\over\log t}\,=\,1/2\qquad
\Q\hbox{-a.s.}
\end{equation}
As a consequence,
$$\lim_{t\to\infty} \left(\log\left(1/|\Pi_1(t)|\right)-t\Phi'(\bar
p)+a\right)
\exp\left(t\Phi(\bar p)\right)|\Pi_1(t)|^{\bar p+1}\,=\,0\qquad
\hbox{$\Q$-a.s.,}$$
and since our goal is to check (\ref{eqliminf}), this enables us to focus
henceforth on
\begin{eqnarray*}
\tilde M(t)\,&=&\,M_a(t)-\left(\log\left(1/|\Pi_1(t)|\right)-t\Phi'(\bar
p)+a\right)
\exp\left(t\Phi(\bar p)\right)|\Pi_1(t)|^{\bar p+1}\\
&=&\,\sum_{i=2}^{\infty}
\left(\log\left(1/|\Pi_i(t)|\right)-t\Phi'(\bar p)+a\right)
\exp\left(t\Phi(\bar p)\right)|\Pi_i(t)|^{\bar p+1}\qquad
\hbox{$\Q$-a.s.,}
\end{eqnarray*}

Next, we compute the conditional expectation of this quantity given
$\g_1(\infty)$, the sigma-field generated by the restriction of the
discrete point measure to the fiber $\R_+\times \p\times\{1\}$, as
follows. By construction of the fragmentation
$\Pi$, each block
$\Pi_i(t)$ for $i\geq2$ got separated from $1$ at some instant $r\in
{\cal
D}_1\cap [0,t]$. More precisely, recall that at such an instant $r$, the
block
$\Pi_1(r-)$ splits into $\pi(r)_{\mid \Pi_1(r-)}$, and that
the block after the split which contains $1$ is $\Pi_1(r)=\pi_1(r)\cap
\Pi_1(r-)$. Thus, there is then some index
$j\geq2$ such that $\Pi_i(t)\subseteq \pi_j(r)\cap \Pi_1(r-)$, where
$\pi_j(r)$ stands for the $j$-th block of the partition $\pi(r)$.
In other words, we may consider the partition of $\{2,\ldots\}$ whose
blocks are of the type
$$B(r,j)\,=\,\left\{i\geq2:\Pi_i(t)\subseteq
\pi_j(r)\cap
\Pi_1(r-)\right\}\,,$$
and then
$\left(\Pi_i(t): i\in B(r,j)\right)$ forms a partition of
$\pi_j(r)\cap \Pi_1(r-)$ which we now analyze.

Lemma \ref{P3'}(ii), standard properties of Poisson random measures,
and the very construction of
$\Pi$ entail that for every $r\in[0,t]$ and $j\geq 2$, conditionally on
$r\in {\cal D}_1$, $\Pi_1(r-)$ and $\pi_j(r)$, the partition
$\left(\Pi_i(t): i\in B(r,j)\right)$
can be given in the form
$\tilde \Pi(t-r)_{\mid \pi_j(r)\cap \Pi_1(r-)}$
where $\tilde \Pi$ is a homogeneous fragmentation distributed as $\Pi$
under
$\P$ and is independent of the sigma-field $\g_1(\infty)$.
Recall from (\ref{defm}) that
$$
\E\left(\sum_{i=1}^{\infty}|\Pi_i(t-r)|^{\bar p+1}\right)
\,=\,\exp(-(t-r)\Phi(\bar p))\,,$$
and (taking the derivative)
$$
\E\left(\sum_{i=1}^{\infty}\log(1/|\Pi_i(t-r)|)|\Pi_i(t-r)|^{\bar
p+1}\right)
\,=\,(t-r)\Phi'(\bar p)\exp(-(t-r)\Phi(\bar p))\,.
$$
The analysis above now entails that
\begin{eqnarray*}
& &\,
\Q\left(\sum_{i=2}^{\infty}
\left(\log\left(1/|\Pi_i(t)|\right)-t\Phi'(\bar p)+a\right)
\exp\left(t\Phi(\bar p)\right)|\Pi_i(t)|^{\bar p+1}
\mid \g_1(\infty)\right)\\
&=&\,\sum_{r\in {\cal D}_1\cap
[0,t]}\sum_{j=2}^{\infty}\left(\log\left(1/|\pi_j(r)\cap\Pi_1(r-)|\right)
-r\Phi'(\bar
p)+a\right)
\exp\left(r\Phi(\bar p)\right)|\pi_j(r)\cap\Pi_1(r-)|^{\bar p+1}
\,.
\end{eqnarray*}
Now observe that there is the identity
$$|\pi_j(r)\cap\Pi_1(r-)|=|\pi_j(r)||\Pi_1(r-)|=|\pi_j(r)|\exp(-\xi(r-))$$
for all
$r\in{\cal D}_1$, $\P$-a.s. and hence also
$\Q$-a.s.
Putting the pieces together, we get
$$\Q(\tilde M(t)\mid \g_1(\infty))\,=\,
\sum_{r\in {\cal D}_1\cap
[0,t]} \exp\left\{-(\bar p+1)(\lambda(r-)-a)\right\}\Sigma(r)\,,$$
where
$$\Sigma(r)\,=\,
\sum_{j=2}^{\infty}
\left(\lambda(r-)-\log|\pi_j(r)|\right)|\pi_j(r)|^{\bar p+1}\,.
$$

As pointed out by Lyons \cite{Lyo}, by the conditional Fatou's theorem,
all that we need is to check that $\lim_{t\to\infty} \Q(\tilde M(t)\mid
\g_1(\infty))<\infty$,  $\Q$-a.s. In this direction, we compute the
$(\Q,\g_1(t))$-predictable compensator corresponding to the point process
$\left\{\Sigma(r), r\in{\cal D}_1\right\}$, and we find
\begin{eqnarray*}
& &\,\lambda(r-)^{-1}\int_{\p}
\mu(d\pi)\left(\lambda(r-)-\log|\pi_1|\right)|\pi_1|^{\bar p}
\left(\sum_{j=2}^{\infty}
\left(\lambda(r-)-\log|\pi_j|\right)|\pi_j|^{\bar p+1}\right)\\
&=&\,
\lambda(r-)^{-1}\int_{\s}
\nu(d{\bf s})\left\{\left(\sum_{j=1}^{\infty}
\left(\lambda(r-)-\log|s_j|\right)|s_j|^{\bar p+1}\right)^2-
\sum_{j=1}^{\infty}
\left(\lambda(r-)-\log|s_j|\right)^2|s_j|^{2\bar p}\right\}\,.
\end{eqnarray*}
Using the fact that $\bar p>0$, it is easily seen that this quantity
can be bounded from above by
$C \left(\lambda(r-)+1+1/\lambda(r-)\right)$ for some constant $C$ that
depends only on $\nu$.
So, all that we need now is to verify that the integral
$$\int_{0}^{\infty}\left(\lambda(r)+1+1/\lambda(r)\right)
\exp\left\{-(\bar p+1)(\lambda(r)-a)\right\}dr$$
converges $\Q$-a.s., which is immediate from (\ref{eq6}).
\QED

Finally, we complete the proof of Proposition \ref{P5}.

{\noindent{\bf Proof of (iii):}\hskip10pt}
To ease the reading, let us
denote
$$Y_i (t) = \exp\left(t\Phi(\bar
p)\right)\left(|\Pi(t)|^{\downarrow}_i\right)^{\bar p+1}\,.$$
We first remark that  for all
$i\in\N$
$\lim
_{t \rightarrow \infty} {Y_i (t)} = 0$,
$\P$-a.s.,
and we deduce from (\ref{defMa}) that $M_a (t) \leq -M'(t) + a M(\bar
p,t)$, for
$t$ large enough. 
Taking the limit as $t\to\infty$, we get
\begin{equation}
\label{mam'}
M_a (\infty) \leq - M'(\infty) \,,\qquad \hbox{$\P$-a.s.,}
\end{equation}
which proves that $-M'(\infty) \geq 0$ $\P$-a.s.

From the fragmentation property at time $1$, we may express $M'(1 +t)$
in the form
$$-M'(1 +t) \,=\, \sum_{i,j} Y_i (1) Y_{i,j} (t) \log \frac{1}{ Y_i (1)
Y_{i,j} (t)}\,,$$
where $\{Y_{i,j}(\cdot), j\in\N\}$ for $i=1,\ldots$ are independent copies
of
$\{Y_j(\cdot), j\in\N\}$, which are also independent of
${\cal G}(1)$. This yields
\begin{equation}
\label{1+t}
-M'(1 +t)
\,=\, \sum_i Y_i (1) \left(-M'_i (t)\right) + \sum_i \left(Y_i (1) \log
\frac{1}{ Y_i (1)}\right) M_i (t)
\end{equation}
where $\{ M_i (\cdot) , i \in \N \}$ (respectively,  $\{ M'_i (\cdot) , i
\in
\N
\}$) are independent copies of $M(\bar p, \cdot)$ (respectively, of
$M'(\cdot)$) and independent of
${\cal G}(1)$.
To get rid of the last infinite random combination of martingales
converging to zero, we first establish the following 
technical result :
\begin{equation}
\label{limsum}
\lim_{t \rightarrow \infty}\, \sum_i \left(Y_i (1) \log  \frac{1}{ Y_i
(1)}\right) M_i (t) \,=\, 0 \qquad \hbox{in probability under $\P$.}
\end{equation}

Indeed, because the
$|\Pi(1)|^{\downarrow}_i$ , $i \in \N$ are ranked in the decreasing
order, and their sum is at most $1$, we have
$|\Pi(1)|^{\downarrow}_i
\leq 1/i$ for every $i$, and thus
$Y_i (1) <1$ for $i >\e^{\Phi'(\bar p)}$. The series $-M'(1) = \sum_i Y_i
(1) |\log {Y_i (1)}|$ is absolutely convergent (and in
$L^1$). Therefore, for every $\epsilon > 0$ there exists $k
>\e^{\Phi'(\bar p)}$ such that
$$\E \left(\sum_{k+1} ^\infty Y_i (1) |\log {Y_i (1)}|\right) \leq
\varepsilon^2\,.$$
Since $\E(M_i (t)) = 1$ for all $i$, the Markov inequality enables us to
write
$$\P\left(\sum_{k+1} ^\infty \left(Y_i (1) |\log {Y_i (1)}|\right)
M_i (t) > \epsilon \right)\leq \epsilon\,.$$
Since the sum of the $k$ remaining terms converges $\P$-a.s. to $0$, the
claim (\ref{limsum}) is proved.

Now we are able to complete the proof of (iii). Assume that $\P(M'(\infty)
= 0) > 0$. From (\ref{1+t}) and (\ref{limsum}) we may write
$$M'(\infty) = Y_1 (1) M'_1 (\infty) + Y_2 (1) M'_2 (\infty) + B$$
where $B = \lim_t \sum_{i=3} ^\infty $ is independent of $(M'_1(\infty)
, M'_2(\infty) )$ conditionally on ${\cal G}_1$. Since $\P(M'(\infty)
\leq 0) = 1$, this entails $\P( B \leq 0) = 1$ and
$$M'(\infty) \leq Y_1 (1) M'_1 (\infty) + Y_2 (1) M'_2 (\infty)\,.$$
This implies $\P (M'(\infty) = 0) \leq \P (M'(\infty) = 0)^2$,
so we would have $\P (M'(\infty) = 0)=1$, which contradicts (ii).
\QED

\vskip 1cm
\noindent {\bf Acknowledgment.} We should like to thank an anonymous referee
for his insightful comments on a first draft of this work.

\end{document}